\documentclass[12pt]{amsart}
\usepackage{amssymb}
\numberwithin{equation}{section}
\begin{document}
\theoremstyle{plain}
\newtheorem{thm}{Theorem}[section]
\newtheorem{theorem}[thm]{Theorem}
\newtheorem{lemma}[thm]{Lemma}
\newtheorem{corollary}[thm]{Corollary}
\newtheorem{proposition}[thm]{Proposition}
\newtheorem{conjecture}[thm]{Conjecture}
\theoremstyle{definition}
\newtheorem{remark}[thm]{Remark}
\newtheorem{remarks}[thm]{Remarks}
\newtheorem{definition}[thm]{Definition}
\newtheorem{example}[thm]{Example}

\newcommand{\sA}{{\mathcal A}}
\newcommand{\sB}{{\mathcal B}}
\newcommand{\sC}{{\mathcal C}}
\newcommand{\sD}{{\mathcal D}}
\newcommand{\sE}{{\mathcal E}}
\newcommand{\sF}{{\mathcal F}}
\newcommand{\sG}{{\mathcal G}}
\newcommand{\sH}{{\mathcal H}}
\newcommand{\sI}{{\mathcal I}}
\newcommand{\sJ}{{\mathcal J}}
\newcommand{\sK}{{\mathcal K}}
\newcommand{\sL}{{\mathcal L}}
\newcommand{\sM}{{\mathcal M}}
\newcommand{\sN}{{\mathcal N}}
\newcommand{\sO}{{\mathcal O}}
\newcommand{\sP}{{\mathcal P}}
\newcommand{\sQ}{{\mathcal Q}}
\newcommand{\sR}{{\mathcal R}}
\newcommand{\sS}{{\mathcal S}}
\newcommand{\sT}{{\mathcal T}}
\newcommand{\sU}{{\mathcal U}}
\newcommand{\sV}{{\mathcal V}}
\newcommand{\sW}{{\mathcal W}}
\newcommand{\sX}{{\mathcal X}}
\newcommand{\sY}{{\mathcal Y}}
\newcommand{\sZ}{{\mathcal Z}}
\newcommand{\ssA}{{\mathfrak A}}
\newcommand{\ssB}{{\mathfrak B}}
\newcommand{\ssC}{{\mathfrak C}}
\newcommand{\ssD}{{\mathfrak D}}
\newcommand{\ssE}{{\mathfrak E}}
\newcommand{\ssF}{{\mathfrak F}}
\newcommand{\ssG}{{\mathfrak G}}
\newcommand{\ssH}{{\mathfrak H}}
\newcommand{\ssI}{{\mathfrak I}}
\newcommand{\ssJ}{{\mathfrak J}}
\newcommand{\ssK}{{\mathfrak K}}
\newcommand{\ssL}{{\mathfrak L}}
\newcommand{\ssM}{{\mathfrak M}}
\newcommand{\ssN}{{\mathfrak N}}
\newcommand{\ssO}{{\mathfrak O}}
\newcommand{\ssP}{{\mathfrak P}}
\newcommand{\ssQ}{{\mathfrak Q}}
\newcommand{\ssR}{{\mathfrak R}}
\newcommand{\ssS}{{\mathfrak S}}
\newcommand{\ssT}{{\mathfrak T}}
\newcommand{\ssU}{{\mathfrak U}}
\newcommand{\ssV}{{\mathfrak V}}
\newcommand{\ssW}{{\mathfrak W}}
\newcommand{\ssX}{{\mathfrak X}}
\newcommand{\ssY}{{\mathfrak Y}}
\newcommand{\ssZ}{{\mathfrak Z}}

\newcommand{\A}{{\mathbb A}}
\newcommand{\B}{{\mathbb B}}
\newcommand{\C}{{\mathbb C}}
\newcommand{\D}{{\mathbb D}}
\newcommand{\E}{{\mathbb E}}
\newcommand{\F}{{\mathbb F}}
\newcommand{\G}{{\mathbb G}}
\newcommand{\HH}{{\mathbb H}}
\newcommand{\I}{{\mathbb I}}
\newcommand{\J}{{\mathbb J}}
\newcommand{\K}{{\mathbb K}}
\newcommand{\LL}{{\mathbb L}}
\newcommand{\M}{{\mathbb M}}
\newcommand{\N}{{\mathbb N}}
\newcommand{\OO}{{\mathbb O}}
\newcommand{\PP}{{\mathbb P}}
\newcommand{\Q}{{\mathbb Q}}
\newcommand{\R}{{\mathbb R}}
\newcommand{\T}{{\mathbb T}}
\newcommand{\U}{{\mathbb U}}
\newcommand{\V}{{\mathbb V}}
\newcommand{\W}{{\mathbb W}}
\newcommand{\X}{{\mathbb X}}
\newcommand{\Y}{{\mathbb Y}}
\newcommand{\Z}{{\mathbb Z}}
\newcommand{\id}{{\rm id}}
\title[Some Results of the Mari\~{n}o-Vafa Formula]{Some Results of the Mari\~{n}o-Vafa Formula}
\author[Yi Li ]{Yi Li}
\address{Center of Mathematical Sciences \\Zhejiang University \\Hangzhou, 310027, China }
\email{yili@cms.zju.edu.cn}
\date{April 20, 2006.}
\begin{abstract}In this paper we derive some new Hodge integral
identities by taking the limits of Mari\~{n}o-Vafa formula. These
identities include the formula of
$\lambda_{1}\lambda_{g}$-integral on $\overline{\sM}_{g,1}$, the
vanishing result of $\lambda_{g}{\rm ch}_{2l}(\E)$-integral on
$\overline{\sM}_{g,1}$ for $1\leq l\leq g-3$. Using the
differential equation of Hodge integrals, we give a recursion
formula of $\lambda_{g-1}$-integrals. Finally, we give two simple
proofs of $\lambda_{g}$ conjecture and some examples of low genus
integral.
\end{abstract}

\maketitle

\setcounter{tocdepth}{5} \setcounter{page}{1}


\section{Introduction}

Based on string duality, Mari\~{n}o and Vafa [10] conjectured a
closed formula on certain Hodge integrals in terms of
representations of symmetric groups. Recently, C.C. Liu, K. Liu
and J. Zhou [6] proved this formula and derived some consequences
from it [7]. In this paper we follow their method to derive some
new Hodge integral identities. One of the main results of this
paper is the following identity: if $1\leq m\leq2g-3$, then
\begin{eqnarray}
&
&-(2g-2-m)!(-1)^{2g-3-m}\int_{\overline{\sM}_{g,1}}\lambda_{g}{\rm
ch}_{2g-2-m}(\E)\psi^{m}_{1} \\
&=&b_{g}\sum^{m-1}_{k=0}\frac{(-1)^{2g-1-k}\binom
{2g-1}{k}\binom {2g-1-k}{2g-1-m}}{2g-1-k}B_{2g-1-m}\nonumber \\
&+&\frac{1}{2}\sum_{g_{1}+g_{2}=g,g_{1},g_{2}>0}b_{g_{1}}b_{g_{2}}\sum^{{\rm
min}(2g_{2}-1,m-1)}_{k=0}\frac{(-1)^{2g_{2}-1-k}\binom
{2g_{2}-1}{k}\binom {2g-1-k}{2g-1-m}}{2g-1-k}B_{2g-1-m}.\nonumber
\end{eqnarray}
As a consequence, we find a new Hodge integral identity: if
$g\geq2$, then
\begin{eqnarray}
\int_{\overline{\sM}_{g,1}}\lambda_{1}\lambda_{g}\psi^{2g-3}_{1}=\frac{1}{12}\left[g(2g-3)b_{g}+b_{1}b_{g-1}\right],
\end{eqnarray}
and also a vanishing result: if $g\geq 2$, then for any $1\leq t
\leq g-1$, we have
\begin{equation}
\int_{\overline{\sM}_{g,1}}\lambda_{g}{\rm
ch}_{2t}(\E)\psi^{2(g-1-t)}_{1}=0.
\end{equation}
Recently, Liu-Xu [8] derived a generalized formula for Hodge
integrals of type (2) by using the $\lambda_g$ conjecture.

The rest of this paper is organized as follows: In Section 2, we
recall the Mari\~{n}o-Vafa formula and the Mumford's relations. In
Section 3, we prove our main theorem and derive a new Hodge
integral identity. In Section 4, we give another simple proof of
$\lambda_{g}$ conjecture. In Section 5, we derive a recursion
formula of $\lambda_{g-1}$-integrals. In the last section, we list
some low genus examples.

\section{Preliminaries}

\subsection{Partitions}A partition $\mu$ of a positive integer $d$ is a
sequence of integers
$\mu_{1}\geq\mu_{2}\geq\cdots\geq\mu_{l(\mu)}>0$ such that
$$\mu_{1}+\cdots+\mu_{l(\mu)}=d=|\mu|,$$
for each positive integer $i$, let
$$m_{i}(\mu)=|\{j|\mu_{j}=i, 1\leq j\leq l(\mu)\}|.$$
The automorphism group ${\rm Aut}(\mu)$ of $\mu$ consists of
possible permutations among the $\mu_{i}$'s, hence its order is
given by
$$|{\rm Aut}(\mu)|=\prod_{i}m_{i}(\mu)!,$$
define the numbers
$$\kappa_{\mu}=\sum^{l(\mu)}_{i=1}\mu_{i}(\mu_{i}-2i+1), \ \ \ \ \ \
z_{\mu}=\prod_{j}m_{j}(\mu)!j^{m_{j}(\mu)}.$$ The Young diagram of
$\mu$ has $l(\mu)$ rows of adjacent squares: the $i$-th row has
$\mu_{i}$ squares. The diagram of $\mu$ can be defines as the set
of points $(i,j)\in\Z\times\Z$ such that $1\leq j\leq \mu_{i}$,
the conjugate of a partition $\mu$ is the partition $\mu'$ whose
diagram is the transpose of the diagram $\mu$. Finally, we
introduce the hook length of $\mu$ at the squre $x\in(i,j)$:
$$h(x)=\mu_{i}+\mu'_{j}-i-j+1.$$
Each partition $\mu$ of $d$ corresponds to a conjugacy class
$C(\mu)$ of the symmetric group $S_{d}$ and each partition $\nu$
corresponds to an irreducible representation $R_{\nu}$ of $S_{d}$,
let $\chi_{\nu}(C(\mu))=\chi_{R_{\nu}}(C(\mu))$ be the value of
the character $\chi_{R_{\nu}}$ on the conjugacy class $C(\mu)$.

\subsection{Mari\~{n}o-Vafa formula}Let
$\overline{\sM}_{g,n}$ denote the Deligne-Mumford moduli stack of
stable curves of genus $g$ with $n$ marked points. Let $\pi:
\overline{\sM}_{g,n+1}\to \overline{\sM}_{g,n}$ be the universal
curve, and let $\omega_{\pi}$ be the relative dualizing sheaf. The
Hodge bundle
$$\E=\pi_{*}\omega_{\pi}$$
is a rank $g$ vector bundle over $\overline{\sM}_{g,n}$ whose
fiber over $[(C,x_{1},\cdots,x_{n})]\in\overline{\sM}_{g,n}$ is
$H^{0}(C,\omega_{C})$, the complex vector space of holomorphic one
forms on $C$. Let $s_{i}: \overline{\sM}_{g,n}\to
\overline{\sM}_{g,n+1}$ denote the section of $\pi$ which
corresponds to the $i$-th marked point, and let
$$\LL_{i}=s^{*}_{i}\omega_{\pi}$$
be the line bundle over $\overline{\sM}_{g,n}$ whose fiber over
$[(C,x_{1},\cdots,x_{n})]\in\overline{\sM}_{g,n}$ is the cotangent
line $T^{*}_{x_{i}}C$ at the $i$-th marked point $x_{i}$. Consider
the Hodge integral
\begin{equation}
\int_{\overline{\sM}_{g,n}}\psi^{j_{1}}_{1}\cdots\psi^{j_{n}}_{n}\lambda^{k_{1}}_{1}\cdots\lambda^{k_{g}}_{g}
\end{equation}
where $\psi_{i}=c_{1}(\LL_{i})$ is the first Chern class of
$\LL_{i}$, and $\lambda_{j}=c_{j}(\E)$ is the $j$-th Chern class
of $\E$. The dimension of $\overline{\sM}_{g,n}$ is $3g-3+n$,
hence (4) is equal to zero unless
$\sum^{n}_{i=1}j_{i}+\sum^{g}_{i=1}ik_{i}=3g-3+n$. Let
\begin{equation}
\Lambda^{\vee}_{g}(u)=u^{g}-\lambda_{1}u^{g-1}+\cdots+(-1)^{g}\lambda_{g}=\sum^{g}_{i=0}(-1)^{i}\lambda_{i}u^{g-i}
\end{equation}
be the Chern polynomial of the dual bundle $\E^{\vee}$ of $\E$.
For any partition $\mu:
\mu_{1}\geq\mu_{2}\geq\cdots\geq\mu_{l(u)}>0$, define
\begin{eqnarray}
\sC_{g,\mu}(\tau)&=&-\frac{\sqrt{-1}^{|\mu|+l(\mu)}}{|{\rm
Aut}(\mu)|}[\tau(\tau+1)]^{l(\mu)-1}\prod^{l(\mu)}_{i=1}\frac{\prod^{\mu_{i}-1}_{a=1}(\mu_{i}\tau+a)}{(\mu_{i}-1)!}
\\
&\cdot&\int_{\overline{\sM}_{g,n}}\frac{\Lambda^{\vee}_{g}(1)\Lambda^{\vee}_{g}(-\tau-1)\Lambda^{\vee}_{g}(\tau)}{\prod^{l(\mu)}_{i=1}(1-\mu_{i}\psi_{i})},
\nonumber\\
\sC_{\mu}(\lambda;\tau)&=&\sum_{g\geq0}\lambda^{2g-2+l(\mu)}\sC_{g,\mu}(\tau),
\end{eqnarray}
here $\tau$ is a formal variable. Note that
\begin{equation*}
\int_{\overline{\sM}_{g,n}}\frac{\Lambda^{\vee}_{g}(1)\Lambda^{\vee}_{g}(-\tau-1)\Lambda^{\vee}_{g}(\tau)}{\prod^{l(\mu)}_{i=1}(1-\mu_{i}\psi_{i})}=|\mu|^{l(\mu)-3}
\end{equation*}
for $l(\mu)\geq 3$, and we use this expression to extend the
definition to the case $l(\mu)<3$.

Introduce formal variables $p=(p_{1},p_{2},\cdots,p_{n},\cdots)$,
and define
$$p_{\mu}=p_{\mu_{1}}\cdots p_{\mu_{l(\mu)}}$$
for a partition $\mu$. Define generating functions
\begin{eqnarray}
\sC(\lambda;\tau,p)&=&\sum_{|\mu|\geq1}\sC_{\mu}(\lambda;\tau)p_{\mu},
\\
\sC(\lambda;\tau,p)^{\bullet}&=&e^{\sC(\lambda;\tau,p)}.
\end{eqnarray}
In [6], Chiu-chu Melissa Liu, Kefeng Liu and Jian Zhou have proved
the following formula which was conjectured by Mari\~{n}o and Vafa
in [10].

\begin{theorem}{\rm {\bf (Mari\~{n}o-Vafa Formula)}}For every partition
$\mu$, we have
\begin{eqnarray*}
\sC(\lambda;\tau,p)&=&\sum_{n\geq1}\frac{(-1)^{n-1}}{n}\sum_{\mu}\left(\sum_{\cup^{n}_{i=1}\mu^{i}=\mu}\prod^{n}_{i=1}\sum_{|\nu^{i}|=|\mu^{i}|}\frac{\chi_{\nu^{i}}(C(\mu^{i}))}{z_{\mu^{i}}}e^{\sqrt{-1}\left(\tau+\frac{1}{2}\right)\kappa_{\nu^{i}}\lambda/2}V_{\nu^{i}}(\lambda)\right)p_{\mu},
\\
\sC(\lambda;\tau,p)^{\bullet}&=&\sum_{|\mu|\geq0}\left(\sum_{|\nu|=|\mu|}\frac{\chi_{\nu}(C(\mu))}{z_{\mu}}e^{\sqrt{-1}\left(\tau+\frac{1}{2}\right)\kappa_{\nu}\lambda/2}V_{\nu}(\lambda)\right)p_{\mu},
\end{eqnarray*}
where
\begin{equation*}
V_{\nu}(\lambda)=\prod_{1\leq a<b\leq l(\nu)}\frac{{\rm
sin}[(\nu_{a}-\nu_{b}+b-a)\lambda/2]}{{\rm sin}[(b-a)\lambda/2]}
\frac{1}{\prod^{l(\nu)}_{i=1}\prod^{\nu_{i}}_{v=1}2{\rm
sin}[(v-i+l(\nu))\lambda/2]}.
\end{equation*}
\end{theorem}
It is known that
$$V_{\nu}(\lambda)=\frac{1}{2^{l(\nu)}\prod_{x\in\nu}{\rm sin}[h(x)\lambda/2]}.$$

\subsection{Mumford's relations}Let
$c_{t}(\E)=\sum^{g}_{i=0}t^{i}\lambda_{i}$, then we have
$$c_{-t}(\E)=t^{g}\Lambda^{\vee}_{g}\left(\frac{1}{t}\right).$$
Mumford's relations [11] are given by
\begin{equation}
c_{t}(\E)c_{-t}(\E)=1,
\end{equation}
equivalently
\begin{equation}
\Lambda^{\vee}_{g}(t)\Lambda^{\vee}_{g}(-t)=(-1)^{g}t^{2g},
\end{equation}
then
\begin{equation}
\lambda^{2}_{k}=\sum^{k}_{i=1}(-1)^{i+1}2\lambda_{k-i}\lambda_{k+i},
\end{equation}
where $\lambda_{0}=1$ and $\lambda_{k}=0$ for $k>g$. Let
$x_{1},\cdots,x_{g}$ be the formal Chern roots of $\E$, the Chern
character is defined by
$${\rm
ch}(\E)=\sum^{g}_{i=1}e^{x_{i}}=g+\sum^{+\infty}_{n=1}\sum^{g}_{i=1}\frac{x^{n}_{i}}{n!},$$
we write
\begin{eqnarray}
{\rm ch}_{0}(\E)&=&g, \\
{\rm ch}_{n}(\E)&=&\frac{1}{n!}\sum^{g}_{i=1}x^{n}_{i},  \ \ \ \ \
\ n=1,2,\cdots.
\end{eqnarray}
From the above identities we have the relation between ${\rm
ch}_{n}(\E)$ and $\lambda_{n}$:
\begin{eqnarray}
n!{\rm ch}_{n}(\E)&=&\sum_{i+j=n}(-1)^{i-1}i\lambda_{i}\lambda_{j},  \ \ \ \ \ \ n<2g,\\
{\rm ch}_{n}(\E)&=&0, \ \ \ \ \ \ n\geq2g.
\end{eqnarray}
It is easy to see that
\begin{eqnarray*}
(2g-1)!{\rm ch}_{2g-1}(\E)&=&(-1)^{g-1}\lambda_{g-1}\lambda_{g}, \\
(2g-2)!{\rm ch}_{2g-2}(\E)&=&(-1)^{g-1}\left((2g-2)\lambda_{g-2}\lambda_{g}-(g-1)\lambda^{2}_{g-1}\right), \\
(2g-3)!{\rm
ch}_{2g-3}(\E)&=&(-1)^{g-1}\left(3\lambda_{g-3}\lambda_{g}-\lambda_{g-1}\lambda_{g-2}\right).
\end{eqnarray*}

\subsection{Bernoulli numbers}The Bernoulli numbers $B_{m}$ are
defined by the following series expansion:
\begin{eqnarray}
\frac{t}{e^{t}-1}=\sum^{+\infty}_{m=0}B_{m}\frac{t^{m}}{m!},
\end{eqnarray}
the first few terms are given by
$$B_{0}=1, \ \ \ B_{1}=-\frac{1}{2}, \ \ \ B_{2}=\frac{1}{6}, \ \ \
B_{3}=0, \ \ \ B_{4}=-\frac{1}{30}, \ \ \ B_{5}=0, \ \ \
B_{6}=\frac{1}{42}.$$ Finally we recall two formulas which will be
used later:
\begin{eqnarray}
\frac{t/2}{{\rm
sin}(t/2)}&=&1+\sum^{+\infty}_{g=1}\frac{2^{2g-1}-1}{2^{2g-1}}\frac{|B_{2g}|}{(2g)!}t^{2g},
\\
\sum^{d-1}_{i=1}i^{m}&=&\sum^{m}_{k=0}\frac{\binom
{m+1}{k}}{m+1}B_{k}d^{m+1-k},
\end{eqnarray}
where $m$ is a positive integer.

\section{Some New Results from Mari\~{n}o-Vafa Formula}

In this section we derive some new results from the
Mari\~{n}o-Vafa formula, we will need two formulas in [7, 2.1 and
5.1].
\begin{theorem}We have the following results:
\begin{eqnarray}
&
&\sum_{g\geq0}\lambda^{2g}\int_{\overline{\sM}_{g,1}}\frac{\frac{d}{d\tau}|_{\tau=0}\left[\Lambda^{\vee}_{g}(1)\Lambda^{\vee}_{g}(\tau)\Lambda^{\vee}_{g}(-\tau-1)\right]}{1-d\psi_{1}}
\nonumber \\
&=&-\sum^{d-1}_{a=1}\frac{1}{a}\frac{d\lambda/2}{d{\rm
sin}(d\lambda/2)}+\sum_{i+j=d,i,j\neq0}\frac{\lambda^{2}}{8{\rm
sin}(i\lambda/2){\rm sin}(j\lambda/2)},
\end{eqnarray}
\begin{equation}
\frac{d}{d\tau}\Big|_{\tau=0}\left[\Lambda^{\vee}_{g}(1)\Lambda^{\vee}_{g}(\tau)\Lambda^{\vee}_{g}(-\tau-1)\right]=-\lambda_{g-1}-\lambda_{g}\sum_{k\geq0}k!(-1)^{k-1}{\rm
ch}_{k}(\E).
\end{equation}
\end{theorem}

\subsection{The coefficient of $\lambda^{2g}$}Introduce the notation
\begin{equation*}
b_{g}=\left\{
      \begin{array}{cc} 1, & g=0, \\
      \frac{2^{2g-1}-1}{2^{2g-1}}\frac{|B_{2g}|}{(2g)!}, & g>0, \\
      \end{array} %
      \right.
\end{equation*}
then the coefficient of $\lambda^{2g}$ in
$-\sum^{d-1}_{a=1}\frac{1}{a}\frac{d\lambda/2}{d{\rm
sin}(d\lambda/2)}$ is
\begin{equation}
\left(-\sum^{d-1}_{a=1}\frac{1}{a}\right)\cdot b_{g}d^{2g-1}.
\end{equation}

If $g_{1}, g_{2}\geq0$ and $g_{1}+g_{2}=g$, define
\begin{equation}
F_{g_{1},g_{2}}(d)=\sum_{i+j=d,i,j\neq0}i^{2g_{1}-1}j^{2g_{2}-1}.
\end{equation}
In [6] it is showed that if $g_{1}, g_{2}\geq1$, then
\begin{equation}
F_{g_{1},g_{2}}(d)=\sum^{2g_{2}-1}_{k=0}\sum^{2g-2-k}_{l=0}\frac{(-1)^{2g_{2}-1-k}}{2g-1-k}\binom
{2g_{2}-1}{k}\binom {2g-1-k}{l}B_{l}d^{2g-1-l},
\end{equation}
for the rest case we have
\begin{eqnarray}
F_{0,g}(d)&=&\sum_{i+j=d,i,j\neq0}i^{-1}j^{2g-1} \\
&=&\sum^{d-1}_{i=1}i^{-1}(d-i)^{2g-1}\nonumber \\
&=&\sum^{2g-1}_{k=0}(-1)^{2g-1-k}\binom
{2g-1}{k}d^{k}\sum^{d-1}_{i=1}i^{2g-2-k}\nonumber \\
&=&\sum^{2g-3}_{k=0}(-1)^{2g-1-k}\binom
{2g-1}{k}d^{k}\sum^{2g-k-2}_{l=0}\frac{\binom
{2g-1-k}{l}}{2g-1-k}B_{l}d^{2g-1-k-l}\nonumber
\\
&-&(2g-1)d^{2g-2}(d-1)+d^{2g-1}\sum^{d-1}_{i=1}\frac{1}{i}
\nonumber\\
&=&\sum^{2g-2}_{k=0}\sum^{2g-k-2}_{l=0}\binom {2g-1}{k}\binom
{2g-1-k}{l}\frac{(-1)^{2g-1-k}}{2g-1-k}B_{l}d^{2g-1-l}\nonumber
\\
&+&(2g-1)d^{2g-2}+d^{2g-1}\sum^{d-1}_{i=1}\frac{1}{i}.\nonumber
\end{eqnarray}
Note that $F_{0,g}(d)=F_{g,0}(d)$ and
\begin{equation}
\sum_{i+j=d,i,j\neq0}\frac{\lambda^{2}}{8{\rm sin}(i\lambda/2){\rm
sin}(j\lambda/2)}=\frac{1}{2}\sum_{g\geq0}\lambda^{2g}\left(\sum_{g_{1}+g_{2}=g}b_{g_{1}}b_{g_{2}}F_{g_{1},g_{2}}(d)\right).
\end{equation}

\subsection{The Main Theorem}Let
$$\sum_{g\geq0}\lambda^{2g}LHS=\sum_{g\geq0}\lambda^{2g}\int_{\overline{\sM}_{g,1}}\frac{\frac{d}{d\tau}|_{\tau=0}\left[\Lambda^{\vee}_{g}(1)\Lambda^{\vee}_{g}(\tau)\Lambda^{\vee}_{g}(-\tau-1)\right]}{1-d\psi_{1}},$$
and
$$\sum_{g\geq0}\lambda^{2g}RHS=-\sum^{d-1}_{a=1}\frac{1}{a}\frac{d\lambda/2}{d{\rm
sin}(d\lambda/2)}+\sum_{i+j=d,i,j\neq0}\frac{\lambda^{2}}{8{\rm
sin}(i\lambda/2){\rm sin}(j\lambda/2)},$$ then we have
\begin{eqnarray}
LHS&=&-\int_{\overline{\sM}_{g,1}}(\lambda_{g-1}\psi_{1}^{2g-1})d^{2g-1}\\
&-&\sum^{2g-2}_{k=0}\left[(2g-2-k)!(-1)^{2g-3-k}\int_{\overline{\sM}_{g,1}}\lambda_{g}{\rm
ch}_{2g-2-k}(\E)\psi^{k}_{1}\right]d^{k},\nonumber \\
RHS&=&-\sum^{d-1}_{a=1}\frac{b_{g}}{a}d^{2g-1}+b_{g}F_{0,g}(d)+\frac{1}{2}\sum_{g_{1}+g_{2}=g,g_{1},g_{2}>0}b_{g_{1}}b_{g_{2}}F_{g_{1},g_{2}}(d).
\end{eqnarray}

Hence we can derive our main theorem:
\begin{theorem}If $1\leq m\leq2g-3$ and $g\geq2$, then
\begin{eqnarray*}
&
&-(2g-2-m)!(-1)^{2g-3-m}\int_{\overline{\sM}_{g,1}}\lambda_{g}{\rm
ch}_{2g-2-m}(\E)\psi^{m}_{1} \\
&=&b_{g}\sum^{m-1}_{k=0}\frac{(-1)^{2g-1-k}}{2g-1-k}\binom
{2g-1}{k}\binom {2g-1-k}{2g-1-m}B_{2g-1-m} \\
&+&\frac{1}{2}\sum_{g_{1}+g_{2}=g,g_{1},g_{2}>0}b_{g_{1}}b_{g_{2}}\sum^{{\rm
min}(2g_{2}-1,m-1)}_{k=0}\frac{(-1)^{2g_{2}-1-k}}{2g-1-k}\binom
{2g_{2}-1}{k}\binom {2g-1-k}{2g-1-m}B_{2g-1-m}.
\end{eqnarray*}
\end{theorem}

\begin{remark}Liu-Liu-Zhou [7] have only considered the cases $m=2g-1$ and
$m= 1$.
\end{remark}

\subsection{The case of $m=2g-3$}If $m=2g-3$, we find that $1!{\rm
ch}_{1}(\E)=\lambda_{1}$, then
\begin{eqnarray*}
LHS&=&-\int_{\overline{\sM}_{g,1}}\lambda_{g}\lambda_{1}\psi^{2g-3}_{1},
\\
RHS&=&b_{g}\sum^{2g-4}_{k=0}\frac{(-1)^{2g-1-k}}{2g-1-k}\binom
{2g-1}{k}\binom {2g-1-k}{2}B_{2} \\
&+&\frac{1}{2}\sum_{g_{1}+g_{2}=g,g_{1},g_{2}>0}b_{g_{1}}b_{g_{2}}\sum^{{\rm
min}(2g_{2}-1,2g-4)}_{k=0}\frac{(-1)^{2g_{2}-1-k}}{2g-1-k}\binom
{2g_{2}-1}{k}\binom {2g-1-k}{2}B_{2}.
\end{eqnarray*}

From the above formula we obtain a new result of the Hodge
integral.
\begin{theorem}If $g\geq2$, then
\begin{eqnarray}
\int_{\overline{\sM}_{g,1}}\lambda_{1}\lambda_{g}\psi^{2g-3}_{1}=\frac{1}{12}\left[g(2g-3)b_{g}+b_{1}b_{g-1}\right].
\end{eqnarray}
\end{theorem}
\begin{proof}Note that
\begin{eqnarray*}
&
&\int_{\overline{\sM}_{g,1}}\lambda_{1}\lambda_{g}\psi^{2g-3}_{1}=-b_{g}B_{2}\sum^{2g-4}_{k=0}\frac{(-1)^{2g-1-k}}{2g-1-k}\binom
{2g-1}{k}\binom {2g-1-k}{2}\\
&-&\frac{B_{2}}{2}\sum_{g_{1}+g_{2}=g,g_{1},g_{2}>0}b_{g_{1}}b_{g_{2}}\sum^{{\rm
min}(2g_{2}-1,2g-4)}_{k=0}\frac{(-1)^{2g_{2}-1-k}}{2g-1-k}\binom
{2g_{2}-1}{k}\binom {2g-1-k}{2},
\end{eqnarray*}
let us write
\begin{eqnarray*}
A_{1}&=&\sum^{2g-4}_{k=0}\frac{(-1)^{2g-1-k}}{2g-1-k}\binom
{2g-1}{k}\binom {2g-1-k}{2}, \\
f_{1}(x)&=&\sum^{2g-3}_{k=0}(-1)^{2g-1-k}\binom
{2g-1}{k}(2g-2-k)x^{2g-3-k}, \\
g_{1}(x)&=&\sum^{2g-2}_{k=0}(-1)^{2g-1-k}\binom
{2g-1}{k}x^{2g-2-k},
\end{eqnarray*}
then
\begin{equation*}
A_{1}=\frac{1}{2}\sum^{2g-4}_{k=0}(-1)^{2g-1-k}\binom
{2g-1}{k}(2g-2-k), \ xg_{1}(x)=(1-x)^{2g-1}-1, \
f_{1}(x)=g_{1}'(x).
\end{equation*}
Hence
\begin{equation*}
f_{1}(x)=\frac{(2g-1)x(1-x)^{2g-2}-(1-x)^{2g-1}+1}{x^{2}}, \ \ \
f_{1}(1)=1,
\end{equation*}
and we obtain
\begin{equation*}
A_{1}=\frac{1}{2}\left[f_{1}(1)-\binom
{2g-1}{2g-3}\right]=-\frac{1}{2}\left[\binom
{2g-1}{2g-3}-1\right].
\end{equation*}
Similarly, we write
$$A_{2}=\sum^{{\rm
min}(2g_{2}-1,2g-4)}_{k=0}\frac{(-1)^{2g_{2}-1-k}}{2g-1-k}\binom
{2g_{2}-1}{k}\binom {2g-1-k}{2},$$ then
\begin{equation*}
A_{2}=\left\{
          \begin{array}{cc}
          \frac{1}{2}\sum^{2g_{2}-1}_{k=0}(-1)^{2g_{2}-1-k}(2g-2-k)\binom {2g_{2}-1}{k}, & {g_{2}\leq g-2}, \\
          & \\
          \frac{1}{2}\sum^{2g-4}_{k=0}(-1)^{2g-1-k}(2g-2-k)\binom {2g-3}{k}, &
          g_{2}=g-1. \\
          \end{array}
        \right.
\end{equation*}
{\bf Case 1: $g_{2}\geq g-2$.} Let
\begin{eqnarray*}
f_{2}(x)&=&\sum^{2g_{2}-1}_{k=0}(-1)^{2g_{2}-1-k}(2g-2-k)\binom
{2g_{2}-1}{k}x^{2g-3-k}, \\
g_{2}(x)&=&\sum^{2g_{2}-1}_{k=0}(-1)^{2g_{2}-1-k}\binom
{2g_{2}-1}{k}x^{2g-2-k}.
\end{eqnarray*}
Since $g\geq g_{2}+2$, then $2g-3-(2g_{2}-1)\geq2>0$ and
$g'_{2}(x)=f_{2}(x)$. On the other hand
$$g_{2}(x)=(1-x)^{2g_{2}-1}x^{2g_{1}-1},$$
hence
$$g'_{2}(x)=-(2g_{2}-1)(1-x)^{2g_{2}-2}x^{2g_{1}-1}+(2g_{1}-1)(1-x)^{2g_{2}-1}x^{2g_{1}-2}$$
and
\begin{equation*}
f_{2}(1)=\left\{
         \begin{array}{cc} -1, & g_{2}=1, \\
         0, & 1<g_{2}\leq g-2. \\
         \end{array} %
         \right.
\end{equation*}
{\bf Case 2: $g_{2}=g-1$.} let
\begin{eqnarray*}
f_{3}(x)&=&\sum^{2g-4}_{k=0}(-1)^{2g-1-k}(2g-2-k)\binom
{2g-3}{k}x^{2g-3-k},\\
g_{3}(x)&=&\sum^{2g-4}_{k=0}(-1)^{2g-1-k}\binom
{2g-3}{k}x^{2g-2-k}.
\end{eqnarray*}
Since $2g-3-(2g-4)=1>0$,
$$\frac{g_{3}(x)}{x}=\sum^{2g-4}_{k=0}(-1)^{2g-3-k}\binom
{2g-3}{k}x^{2g-3-k}=(1-x)^{2g-3}-1$$ and
$$g'_{3}(x)=(1-x)^{2g-3}-1-(2g-3)(1-x)^{2g-4}x,$$
therefore we have
\begin{equation*}
f_{3}(1)=\left\{
         \begin{array}{cc} -2, & g=2, \\
          -1, & g>2. \\
          \end{array} %
          \right.
\end{equation*}
From the values of $f_{1}(1), f_{2}(1), f_{3}(1)$, we obtain
\begin{eqnarray*}
&
&\int_{\overline{\sM}_{g,1}}\lambda_{1}\lambda_{g}\psi^{2g-3}_{1}
\\
&=&-b_{g}B_{2}A_{1}-\frac{B_{2}}{2}\left[\frac{1}{2}\sum_{g_{1}+g_{2}=g,
1\leq g_{2}\leq
g-2}b_{g_{1}}b_{g_{2}}f_{2}(1)+\frac{1}{2}b_{1}b_{g-1}f_{3}(1)\right]
\\
&=&\frac{B_{2}}{2}[-b_{g}A_{1}+b_{1}b_{g-1}]
\\
&=&\frac{1}{12}\left[g(2g-3)b_{g}+b_{1}b_{g-1}\right]
\end{eqnarray*}
\end{proof}

Since $B_{n}=0$ for $n$ odd and $n>1$, we also have the following
vanishing result.

\begin{theorem}If $g\geq 2$, then for any $1\leq t \leq g-1$, we
have
\begin{equation}
\int_{\overline{\sM}_{g,1}}\lambda_{g}{\rm
ch}_{2t}(\E)\psi^{2(g-1-t)}_{1}=0.
\end{equation}
\end{theorem}

\section{Another Simple Proof of $\lambda_{g}$ Conjecture}Let $|\mu|=d,
l(\mu)=n$, denote by $[\sC_{g,\mu}(\tau)]_{k}$ the coefficient of
$\tau^{k}$ in the polynomial $\sC_{g,\mu}(\tau)$, and let
\begin{eqnarray*}
J^{0}_{g,\mu}(\tau)&=&\sqrt{-1}^{|\mu|-l(\mu)}\sC_{g,\mu}(\tau), \\
J^{1}_{g,\mu}(\tau)&=&\sqrt{-1}^{|\mu|-l(\mu)-1}\left(\sum_{\nu\in
J(\mu)}I_{1}(\nu)\sC_{g,\nu}(\tau)+\sum_{\nu\in
C(\mu)}I_{2}(\nu)\sC_{g-1,\nu}(\tau)\right. \\
&+&\left.\sum_{g_{1}+g_{2}=g, \nu^{1}\cup\nu^{2}\in
C(\mu)}I_{3}(\nu^{1},\nu^{2})\sC_{g_{1},\nu^{1}}(\tau)\sC_{g_{2},\nu^{2}}(\tau)\right).
\end{eqnarray*}
The set $J(\mu)$ consists of partitions of $d$ of the form
\begin{equation*}
\nu=(\mu_{1},\cdots,\widehat{\mu}_{i},\cdots,\mu_{l(\mu)},\mu_{i}+\mu_{j})
\end{equation*}
and the set $C(\mu)$ consists of partitions of $d$ of the form
\begin{equation*}
\nu=(\mu_{1},\cdots,\widehat{\mu}_{i},\cdots,\mu_{l(\mu)},j,k)
\end{equation*}
where $j+k=\mu_{i}$. The definitions of $I_{1}, I_{2}$ and $I_{3}$
can be found in [5]. Liu-Liu-Zhou [6] have proved the following
differential equation:
\begin{equation}
\frac{d}{d\tau}J^{0}_{g,\mu}(\tau)=-J^{1}_{g,\mu}(\tau).
\end{equation}
It is straightforward to check that
\begin{eqnarray*}
\left[\sC_{g,\mu}(\tau)\right]_{n-1}&=&-\frac{\sqrt{-1}^{d+n}}{|{\rm
Aut}(\mu)|}\int_{\overline{\sM}_{g,n}}\frac{\lambda_{g}}{\prod^{n}_{i=1}(1-\mu_{i}\psi_{i})},
\\
\left[\sum_{\nu\in
J(\mu)}I_{1}(\nu)\sC_{g,\nu}(\tau)\right]_{n-2}&=&-\frac{\sqrt{-1}^{d+n-1}}{|{\rm
Aut(\nu)|}}\int_{\overline{\sM}_{g,n-1}}\frac{\lambda_{g}}{\prod^{n-1}_{i=1}(1-\mu_{i}\psi_{i})},
\end{eqnarray*}
\begin{eqnarray*}
\left[\sum_{\nu\in
C(\mu)}I_{2}(\nu)\sC_{g-1,\nu}(\tau)\right]_{n-2}&=&0, \\
\left[\sum_{g_{1}+g_{2}=g, \nu^{1}\cup\nu^{2}\in
C(\mu)}I_{3}(\nu^{1},\nu^{2})\sC_{g_{1},\nu^{1}}(\tau)\sC_{g_{2},\nu^{2}}(\tau)\right]_{n-2}&=&0,
\end{eqnarray*}
hence, from (31) we have the identity
\begin{equation}
\frac{n-1}{|\rm
Aut(\mu)|}\int_{\overline{\sM}_{g,n}}\frac{\lambda_{g}}{\prod^{n}_{i=1}(1-\mu_{i}\psi_{i})}=\sum_{\nu\in
J(\mu)}\frac{I_{1}(\nu)}{|{\rm
Aut}(\nu)|}\int_{\overline{\sM}_{g,n-1}}\frac{\lambda_{g}}{\prod^{n-1}_{i=1}(1-\nu_{i}\psi_{i})}.
\end{equation}

\begin{theorem}For any partition $\mu:
\mu_{1}\geq\mu_{2}\geq\cdots\geq\mu_{n}>0$ of $d$ and $g>0$, then
\begin{equation}
\int_{\overline{\sM}_{g,n}}\frac{\lambda_{g}}{\prod^{n}_{i=1}(1-\mu_{i}\psi_{i})}=d^{2g+n-3}b_{g}.
\end{equation}
\end{theorem}
\begin{proof}Recall the definition of $I_{1}(\nu)$, where
$\nu=(\mu_{1},\cdots,\widehat{\mu_{i}},\cdots,\widehat{\mu}_{j},\cdots,\mu_{n},\mu_{i}+\mu_{j})$:
$$I_{1}(\nu)=\frac{\mu_{i}+\mu_{j}}{1+\delta^{\mu_{i}}_{\mu_{j}}}m_{\mu_{i}+\mu_{j}}(\nu),$$
and it is easy to see that
$$\frac{m_{\mu_{i}+\mu_{j}}(\nu)}{|{\rm
Aut}(\nu)|}=\frac{m_{\mu_{i}}(\mu)(m_{\mu_{j}}(\mu)-\delta^{\mu_{i}}_{\mu_{j}})}{|{\rm
Aut}(\mu)|}.$$ Let
\begin{equation*}
\mu:
\underbrace{\mu_{k_{1}}=\cdots=\mu_{k_{1}}}_{t_{1}}>\underbrace{\mu_{k_{2}}=\cdots=\mu_{k_{2}}}_{t_{2}}>\cdots>\underbrace{\mu_{k_{s}}=\cdots=\mu_{k_{s}}}_{t_{s}}>0,
\end{equation*}
where
$$\sum^{s}_{i=1}t_{i}=n, \ \ \ \ \ \
\sum^{s}_{i=1}t_{i}\mu_{k_{i}}=d,$$ then
\begin{eqnarray*}
& &\sum_{\nu\in J(\mu)}\frac{I_{1}(\nu)}{|{\rm Aut}(\nu)|}
\\
&=&\frac{1}{|{\rm Aut}(\mu)|}\sum_{\nu\in
J(\mu)}\frac{\mu_{i}+\mu_{j}}{1+\delta^{\mu_{i}}_{\mu_{j}}}[m_{\mu_{i}}(\mu)(m_{\mu_{j}}(\mu)-\delta^{\mu_{i}}_{\mu_{j}})]
\\
&=&\frac{1}{|{\rm
Aut}(\mu)|}\left[\frac{1}{2}\sum^{s}_{i=1}\sum_{j\neq
i}(\mu_{k_{i}}+\mu_{k_{j}})m_{\mu_{k_{i}}}(\mu)m_{\mu_{k_{j}}}(\mu)+\sum^{s}_{i=1}\mu_{k_{i}}m_{\mu_{k_{i}}}(\mu)(m_{\mu_{k_{i}}}(\mu)-1)\right]
\\
&=&\frac{1}{|{\rm
Aut}(\mu)|}\left[\frac{1}{2}\sum^{s}_{i=1}\sum_{j\neq
i}(\mu_{k_{i}}+\mu_{k_{j}})t_{i}t_{j}+\sum^{s}_{i=1}\mu_{k_{i}}t_{i}(t_{i}-1)\right]
\\
&=&\frac{1}{|{\rm Aut}(\mu)|}\left[\sum^{s}_{i=1}\sum_{j\neq
i}\mu_{k_{j}}t_{j}t_{i}+\sum^{s}_{i=1}\mu_{k_{i}}t^{2}_{i}-d\right]
\\
&=&\frac{1}{|{\rm
Aut}(\mu)|}\left[\sum^{s}_{i=1}t_{i}(d-\mu_{k_{i}}t_{i})+\sum^{s}_{i=1}\mu_{k_{i}}t^{2}_{i}-d\right]
\\
&=&\frac{(n-1)d}{|{\rm Aut}(\mu)|}.
\end{eqnarray*}
By the induction of $n$ and the initial value of the
Mari\~{n}o-Vafa formula
\begin{equation*}
\int_{\overline{\sM}_{g,1}}\frac{\lambda_{g}}{1-\mu_{1}\psi_{1}}=d^{2g-2}b_{g},
\end{equation*}
we have
\begin{equation*}
\int_{\overline{\sM}_{g,n}}\frac{\lambda_{g}}{\prod^{n}_{i=1}(1-\mu_{i}\psi_{i})}=d\cdot
d^{2g+n-1-3}b_{g}=d^{2g+n-3}b_{g}.
\end{equation*}
\end{proof}

\begin{corollary}The following  $\lambda_{g}$ conjecture is true:
\begin{equation}
\int_{\overline{\sM}_{g,n}}\lambda_{g}\prod^{n}_{l=1}\psi^{k_{l}}_{l}=\binom{2g+n-3}{k_{1},\cdots,k_{n}}b_{g},
\end{equation}
where $g>0$.
\end{corollary}

\section{A recursion Formula of the $\lambda_{g-1}$ Integral}

E.Getzler, A.Okounkov and R.Pandharipande have derived explicit
formula for the multipoint series of $\C\PP^{1}$ in degree $0$
from the Toda hierarchy [2], then they obtained certain formulas
for the Hodge integrals
$\int_{\overline{\sM}_{g,n}}\lambda_{g-1}\psi^{k_{1}}_{1}\cdots\psi^{k_{n}}_{n}$.
In this section we give an effective recursion formula of the
$\lambda_{g-1}$ integrals using Mari\~{n}o-Vafa formula. It is
straightforward to check the following lemma.
\begin{lemma}We have the following identities
\begin{eqnarray*}
\left[\prod^{n}_{i=1}\frac{\prod^{\mu_{i}-1}_{a=1}(\mu_{i}\tau+a)}{(\mu_{i}-1)!}\right]_{0}&=&1,
\\
\left[\prod^{n}_{i=1}\frac{\prod^{\mu_{i}-1}_{a=1}(\mu_{i}\tau+a)}{(\mu_{i}-1)!}\right]_{1}&=&\sum^{n}_{i=1}\sum^{\mu_{i}-1}_{a=1}\frac{\mu_{i}}{a},
\\
\left[\Lambda^{\vee}_{g}(1)\Lambda^{\vee}_{g}(-\tau-1)\Lambda^{\vee}_{g}(\tau)\right]_{0}&=&\lambda_{g},
\\
\left[\Lambda^{\vee}_{g}(1)\Lambda^{\vee}_{g}(-\tau-1)\Lambda^{\vee}_{g}(\tau)\right]_{0}&=&-\lambda_{g-1}-\sum_{k\geq0}k!(-1)^{k-1}{\rm
ch}_{k}(\E)\lambda_{g},
\end{eqnarray*}
and
\begin{eqnarray*}
\left[\sC_{g,\mu}(\tau)\right]_{n-1}&=&-\frac{\sqrt{-1}^{d+n}}{|{\rm
Aut}(\mu)|}\int_{\overline{\sM}_{g,n}}\frac{\lambda_{g}}{\prod^{n}_{i=1}(1-\mu_{i}\psi_{i})},
\\
\left[\sC_{g,\mu}(\tau)\right]_{n}&=&-\frac{\sqrt{-1}^{d+n}}{|{\rm
Aut}(\mu)|}\left[n-1+\sum^{n}_{i=1}\sum^{\mu_{i}-1}_{a=1}\frac{\mu_{i}}{a}\right]\int_{\overline{\sM}_{g,n}}\frac{\lambda_{g}}{\prod^{n}_{i=1}(1-\mu_{i}\psi_{i})}
\\
&+&\frac{\sqrt{-1}^{d+n}}{|{\rm
Aut}(\mu)|}\int_{\overline{\sM}_{g,n}}\frac{\lambda_{g-1}+\sum_{k\geq0}k!(-1)^{k-1}{\rm
ch}_{k}(\E)\lambda_{g}}{\prod^{n}_{i=1}(1-\mu_{i}\psi_{i})}.
\end{eqnarray*}
\end{lemma}

Now, we can state our main theorem in this section using equation
(31).
\begin{theorem}For any partition $\mu$ with $l(\mu)=n$, we have the
following recursion formula
\begin{eqnarray*}
& &\frac{n}{|{\rm
Aut}(\mu)|}\left[n-1+\sum^{n}_{i=1}\sum^{\mu_{i}-1}_{a=1}\frac{\mu_{i}}{a}\right]\int_{\overline{\sM}_{g,n}}\frac{\lambda_{g}}{\prod^{n}_{i=1}(1-\mu_{i}\psi_{i})}
\\
&-&\frac{n}{|{\rm
Aut}(\mu)|}\int_{\overline{\sM}_{g,n}}\frac{\lambda_{g-1}+\sum_{k\geq0}k!(-1)^{k-1}{\rm
ch}_{k}(\E)\lambda_{g}}{\prod^{n}_{i=1}(1-\mu_{i}\psi_{i})} \\
&=&\sum_{\nu\in J(\mu)}\frac{I_{1}(\nu)}{|{\rm
Aut}(\nu)|}\left[n-2+\sum^{n-1}_{i=1}\sum^{\nu_{i}-1}_{a=1}\frac{\nu_{i}}{a}\right]\int_{\overline{\sM}_{g,n-1}}\frac{\lambda_{g}}{\prod^{n-1}_{i=1}(1-\nu_{i}\psi_{i})}
\\
&-&\sum_{\nu\in J(\mu)}\frac{I_{1}(\nu)}{|{\rm
Aut}(\nu)|}\int_{\overline{\sM}_{g,n-1}}\frac{\lambda_{g-1}+\sum_{k\geq0}k!(-1)^{k-1}{\rm
ch}_{k}(\E)\lambda_{g}}{\prod^{n-1}_{i=1}(1-\nu_{i}\psi_{i})} \\
&+&\sum_{g_{1}+g_{2}=g,
g_{1},g_{2}\geq0}\sum_{\nu^{1}\cup\nu^{2}\in
C(\mu)}\frac{I_{3}(\nu^{1},\nu^{2})}{|{\rm Aut}(\nu^{1})||{\rm
Aut}(\nu^{2})|}
\\
&\cdot&\int_{\overline{\sM}_{g_{1},n_{1}}}\frac{\lambda_{g_{1}}}{\prod^{n_{1}}_{i=1}(1-\nu^{1}_{i}\psi_{i})}\int_{\overline{\sM}_{g_{2},n_{2}}}\frac{\lambda_{g_{2}}}{\prod^{n_{2}}_{i=1}(1-\nu^{2}_{i}\psi_{i})}.
\end{eqnarray*}
\end{theorem}

\subsection{The $\lambda_{g}$-Integral}In this subsection, we re-derive the $\lambda_{g}$-integral from theorem 5.2.. Let $\mu_{i}=Nx_{i}$ for some $N\in\N$ and
$x_{i}\in\R$, from Kim-Liu[4]'s method and consider the
coefficients of ${\rm ln}NN^{2g+n-2}$ in theorem 5.2., then
\begin{eqnarray*}
&
&n(x_{1}+\cdots+x_{n})\prod^{n}_{l=1}x^{k_{l}}_{l}\int_{\overline{\sM}_{g,n}}\lambda_{g}\prod^{n}_{l=1}\psi^{k_{l}}_{l}
\\
&=&\frac{1}{2}\sum^{n}_{i=1}\sum_{j\neq
i}(x_{i}+x_{j})^{k_{i}+k_{j}}(x_{1}+\cdots+x_{n})\prod_{l\neq
i,j}x^{k_{l}}_{l}\int_{\overline{\sM}_{g,n-1}}\lambda_{g}\psi^{k_{i}+k_{j}-1}\prod_{l\neq
i,j}\psi^{k_{l}}_{l} \\
&+&(x_{1}+\cdots+x_{n})\prod^{n}_{l=1}x^{k_{l}}_{l}\int_{\overline{\sM}_{g,n}}\lambda_{g}\prod^{n}_{l=1}\psi^{k_{l}}_{l},
\end{eqnarray*}
i.e.
\begin{eqnarray}
&
&(n-1)\prod^{n}_{l=1}x^{k_{l}}_{l}\int_{\overline{\sM}_{g,n}}\lambda_{g}\prod^{n}_{l=1}\psi^{k_{l}}_{l}
\\
&=&\frac{1}{2}\sum^{n}_{i=1}\sum_{j\neq
i}(x_{i}+x_{j})^{k_{i}+k_{j}}\prod_{l\neq
i,j}x^{k_{l}}_{l}\int_{\overline{\sM}_{g,n-1}}\lambda_{g}\psi^{k_{i}+k_{j}-1}\prod_{l\neq
i,j}\psi^{k_{l}}_{l}.\nonumber
\end{eqnarray}
After introducing the formal variables $s_{l}\in\R^{+}$ and
applying the Laplace transformation
\begin{equation*}
\int^{+\infty}_{0}x^{k}e^{-x/2s}dx=k!(2s)^{k+1}, \ \ \ s>0,
\end{equation*}
we select the coefficient of $\prod^{n}_{l=1}(2s_{l})^{k_{l}+1}$
from the transformation of (35), then we derive
\begin{equation}
(n-1)\int_{\overline{\sM}_{g,n}}\lambda_{g}\prod^{n}_{l=1}\psi^{k_{l}}_{l}=\frac{1}{2}\sum^{n}_{i=1}\sum_{j\neq
i}\frac{(k_{i}+k_{j})!}{k_{i}!k_{j}!}\int_{\overline{\sM}_{g,n-1}}\lambda_{g}\psi^{k_{i}+k_{j}-1}\prod_{l\neq
i,j}\psi^{k_{l}}_{l}.
\end{equation}
By the induction of $n$, we obtain the $\lambda_{g}$ conjecture
\begin{equation*}
\int_{\overline{\sM}_{g,n}}\lambda_{g}\prod^{n}_{l=1}\psi^{k_{l}}_{l}=\binom{2g+n-3}{k_{1},\cdots,k_{n}}b_{g},
\end{equation*}
in fact, in (36) we have
\begin{eqnarray*}
RHS&=&\frac{1}{2}\sum^{n}_{i=1}\sum_{j\neq
i}\frac{(k_{i}+k_{j})!}{k_{i}!k_{j}!}\frac{(2g+n-4)!}{\prod_{l\neq
i,j}k_{l}!(k_{i}+k_{j}-1)!}b_{g} \\
&=&\frac{1}{2}\sum^{n}_{i=1}\sum_{j\neq
i}\frac{k_{i}+k_{j}}{2g+n-3}\binom{2g+n-3}{k_{1},\cdots,k_{n}}b_{g},
\end{eqnarray*}
note that $k_{1}+\cdots+k_{n}=2g+n-3$, therefore
\begin{eqnarray*}
\frac{1}{2}\sum^{n}_{i=1}\sum_{j\neq
i}(k_{i}+k_{j})&=&\frac{1}{2}\sum^{n}_{i=1}[(n-1)k_{i}+(2g+n-3-k_{i})]
\\
&=&\frac{1}{2}\left[(n-2)(2g+n-3)+(2g+n-3)n\right] \\
&=&\frac{1}{2}[(2n-2)(2g+n-3)] \\
&=&(n-1)(2g+n-3).
\end{eqnarray*}

\subsection{The Recursion Formula of $\lambda_{g-1}$-integral}We have found the {\it singular part} $\sum^{n}_{i=1}\sum^{\mu_{i}-1}_{a=1}\frac{\mu_{i}}{a}$ in theorem 5.2., using
the following theorem, we can eliminate this part and derive the
recursion formula of $\lambda_{g-1}$-integral. The notation
$[F]_{sing}$ means the singular part of $F$. First, in theorem
5.2., we have
\begin{eqnarray*}
\left[\frac{LHS}{d^{2g+n-4}b_{g}}\right]_{sing}&=&n\sum^{n}_{i=1}\sum^{\mu_{i}-1}_{a=1}\frac{\mu_{i}}{a}d,
\\
\left[\frac{RHS}{d^{2g+n-4}b_{g}}\right]_{sing}&=&\frac{1}{2}\sum^{n}_{i=1}\sum_{j\neq
i}(\mu_{i}+\mu_{j})\left[\sum_{l\neq
i,j}\sum^{\mu_{l}-1}_{a=1}\frac{\mu_{l}}{a}+(\mu_{i}+\mu_{j})\sum^{\mu_{i}+\mu_{j}-1}_{a=1}\frac{1}{a}\right]
\\
&+&\sum^{n}_{i=1}\sum^{\mu_{i}-1}_{a=1}\frac{\mu_{i}}{a}d-\sum^{n}_{i=1}\sum_{j\neq
i}\sum^{\mu_{i}+\mu_{j}-1}_{a=\mu_{j}+1}\frac{\mu_{j}(\mu_{i}+\mu_{j})}{a}.
\end{eqnarray*}

\begin{theorem}Under the above notation, we have
\begin{equation*}
\left[\frac{RHS}{d^{2g+n-4}b_{g}}\right]_{sing}=\left[\frac{LHS}{d^{2g+n-4}b_{g}}\right]_{sing}+2(n-1)d.
\end{equation*}
\end{theorem}
\begin{proof}Since
\begin{eqnarray*}
& &\sum^{n}_{i=1}\sum_{j\neq
i}\sum^{\mu_{i}+\mu_{j}-1}_{a=\mu_{j}+1}\frac{\mu_{j}(\mu_{i}+\mu_{j})}{a}
\\
&=&\sum^{n}_{i=1}\sum_{j\neq
i}\sum^{\mu_{i}+\mu_{j}-1}_{a=\mu_{j}+1}\frac{(\mu_{i}+\mu_{j})^{2}}{a}-\sum^{n}_{i=1}\sum_{j\neq
i}\sum^{\mu_{i}+\mu_{j}-1}_{a=\mu_{j}+1}\frac{\mu_{i}(\mu_{i}+\mu_{j})}{a}
\\
&=&\sum^{n}_{i=1}\sum_{j\neq
i}\sum^{\mu_{i}+\mu_{j}-1}_{a=1}\frac{(\mu_{i}+\mu_{j})^{2}}{a}-\sum^{n}_{i=1}\sum_{j\neq
i}\sum^{\mu_{j}}_{a=1}\frac{(\mu_{i}+\mu_{j})^{2}}{a} \\
&-&\sum^{n}_{i=1}\mu_{i}\sum_{j\neq
i}\sum^{\mu_{i}+\mu_{j}-1}_{a=1}\frac{\mu_{i}+\mu_{j}}{a}+\sum^{n}_{i=1}\mu_{i}\sum_{j\neq
i}\sum^{\mu_{j}}_{a=1}\frac{(\mu_{i}+\mu_{j})}{a} \\
&=&\frac{1}{2}\sum^{n}_{i=1}\sum_{j\neq
i}\sum^{\mu_{i}+\mu_{j}-1}_{a=1}\frac{(\mu_{i}+\mu_{j})^{2}}{a}-\sum^{n}_{i=1}\sum_{j\neq
i}\sum^{\mu_{j}-1}_{a=1}\frac{(\mu_{i}+\mu_{j})^{2}}{a} \\
&+&\sum^{n}_{i=1}\mu_{i}\sum_{j\neq
i}\sum^{\mu_{j}-1}_{a=1}\frac{(\mu_{i}+\mu_{j})}{a}-\sum^{n}_{i=1}\sum_{j\neq
i}(\mu_{i}+\mu_{j}),
\end{eqnarray*}
where we use the identity
\begin{eqnarray*}
\sum^{n}_{i=1}\sum_{j\neq
i}\sum^{\mu_{i}+\mu_{j}-1}_{a=1}\frac{(\mu_{i}+\mu_{j})\mu_{i}}{a}&=&\sum^{n}_{i=1}\sum_{j\neq
i}\sum^{\mu_{i}+\mu_{j}-1}_{a=1}\frac{(\mu_{i}+\mu_{j})\mu_{j}}{a}
\\
&=&\frac{1}{2}\sum^{n}_{i=1}\sum_{j\neq
i}\sum^{\mu_{i}+\mu_{j}-1}_{a=1}\frac{(\mu_{i}+\mu_{j})^{2}}{a}.
\end{eqnarray*}
Note that $\sum^{n}_{i=1}\sum_{j\neq i}(\mu_{i}+\mu_{j})=2(n-1)d$,
hence
\begin{eqnarray*}
& &\left[\frac{RHS}{d^{2g+n-4}b_{g}}\right]_{sing}
\\
&=&\frac{1}{2}\sum^{n}_{i=1}\sum_{j\neq
i}(\mu_{i}+\mu_{j})\sum_{l\neq
i,j}\sum^{\mu_{l}-1}_{a=1}\frac{\mu_{l}}{a}+\sum^{n}_{i=1}\sum^{\mu_{i}-1}_{a=1}\frac{\mu_{i}}{a}d
\\
&+&\sum^{n}_{i=1}\sum_{j\neq
i}\sum^{\mu_{j}-1}_{a=1}\frac{(\mu_{i}+\mu_{j})^{2}}{a}-\sum^{n}_{i=1}\sum_{j\neq
i}\mu_{i}\sum^{\mu_{j}-1}_{a=1}\frac{(\mu_{i}+\mu_{j})}{a}+\sum^{n}_{i=1}\sum_{j\neq i}(\mu_{i}+\mu_{j}) \\
&=&\left(\sum^{n}_{i=1}\sum^{\mu_{i}-1}_{a=1}\frac{\mu_{i}}{a}\right)d+\frac{1}{2}\sum^{n}_{i=1}\sum_{j\neq
i}(\mu_{i}+\mu_{j})\sum_{l\neq
i,j}\sum^{\mu_{l}-1}_{a=1}\frac{\mu_{l}}{a}
\\
&+&\sum^{n}_{i=1}\sum_{j\neq
i}\sum^{\mu_{j}-1}_{a=1}\frac{\mu_{j}(\mu_{i}+\mu_{j})}{a}+2(n-1)d,
\end{eqnarray*}
it is straightforward to check that
\begin{eqnarray*}
\frac{1}{2}\sum^{n}_{i=1}\sum_{j\neq
i}(\mu_{i}+\mu_{j})\sum_{l\neq
i,j}\sum^{\mu_{l}-1}_{a=1}\frac{\mu_{l}}{a}&=&(n-2)\sum^{n}_{i=1}\mu_{i}\sum_{j\neq
i}\sum^{\mu_{j}-1}_{a=1}\frac{\mu_{j}}{a}, \\
\sum^{n}_{i=1}\sum_{j\neq
i}\sum^{\mu_{j}-1}_{a=1}\frac{\mu_{j}(\mu_{i}+\mu_{j})}{a}&=&\sum^{n}_{i=1}\mu_{i}\sum_{j\neq
i}\sum^{\mu_{j}-1}_{a=1}\frac{\mu_{j}}{a}+\sum^{n}_{i=1}\sum_{j\neq
i}\mu_{j}\sum^{\mu_{j}-1}_{a=1}\frac{\mu_{j}}{a}, \\
\sum^{n}_{i=1}\sum_{j\neq
i}\mu_{j}\sum^{\mu_{j}-1}_{a=1}\frac{\mu_{j}}{a}&=&(n-1)\sum^{n}_{i=1}\mu_{i}\sum^{\mu_{i}-1}_{a=1}\frac{\mu_{i}}{a}.
\end{eqnarray*}
Finally, we obtain
\begin{eqnarray*}
& &\left[\frac{RHS}{d^{2g+n-4}b_{g}}\right]_{sing}
\\
&=&\left(\sum^{n}_{i=1}\sum^{\mu_{i}-1}_{a=1}\frac{\mu_{i}}{a}\right)d+(n-2)\sum^{n}_{i=1}\mu_{i}\sum_{j\neq
i}\sum^{\mu_{j}-1}_{a=1}\frac{\mu_{j}}{a} \\
&+&\sum^{n}_{i=1}\mu_{i}\sum_{j\neq
i}\sum^{\mu_{j}-1}_{a=1}\frac{\mu_{j}}{a}+(n-1)\sum^{n}_{i=1}\mu_{i}\sum^{\mu_{i}-1}_{a=1}\frac{\mu_{i}}{a}
+2(n-1)d\\
&=&\left(\sum^{n}_{i=1}\sum^{\mu_{i}-1}_{a=1}\frac{\mu_{i}}{a}\right)d+(n-1)\sum^{n}_{i=1}\mu_{i}\left(\sum_{j\neq
i}\sum^{\mu_{i}-1}_{a=1}\frac{\mu_{j}}{a}+\sum^{\mu_{i}-1}_{a=1}\frac{\mu_{i}}{a}\right)+2(n-1)d
\\
&=&\left(\sum^{n}_{i=1}\sum^{\mu_{i}-1}_{a=1}\frac{\mu_{i}}{a}\right)d+(n-1)\left(\sum^{n}_{i=1}\sum^{\mu_{i}-1}_{a=1}\frac{\mu_{i}}{a}\right)d+2(n-1)d
\\
&=&\left[\frac{LHS}{d^{2g+n-4}b_{g}}\right]_{sing}+2(n-1)d.
\end{eqnarray*}
\end{proof}

Let $\R^{k}[\mu_{1},\cdots,\mu_{n}]$ be the space of all
homogeneous polynomials with real coefficients in
$\mu_{1},\cdots,\mu_{n}$ of degree $k$, then it is the subring of
$\R[\mu_{1},\cdots,\mu_{n}]$. From the Theorem 5.3, we obtain the
recursion formula of $\lambda_{g-1}$ Hodge integral.
\begin{theorem}For any partition $\mu$ with $l(\mu)=n$ and $|\mu|=d$, we
have the recursion formula
\begin{eqnarray*}
& &\frac{n}{|{\rm
Aut}(\mu)|}\int_{\overline{\sM}_{g,n}}\frac{\lambda_{g-1}}{\prod^{n}_{i=1}(1-\mu_{i}\psi_{i})}
\\
&=&\sum_{\nu\in J(\mu)}\frac{I_{1}(\nu)}{|{\rm
Aut}(\nu)|}\int_{\overline{\sM}_{g,n-1}}\frac{\lambda_{g-1}}{\prod^{n-1}_{i=1}(1-\nu_{i}\psi_{i})}
\\
&-&\sum_{g_{1}+g_{2}=g, \nu^{1}\cup\nu^{2}\in
C(\mu)}\frac{I_{3}(\nu^{1},\nu^{2})}{|{\rm Aut}(\nu^{1})||{\rm
Aut}(\nu^{2})|}d^{2g_{1}+n_{1}-3}_{1}d^{2g_{2}+n_{2}-3}_{2}b_{g_{1}}b_{g_{2}}.
\end{eqnarray*}
under the ring $\R^{2g-2+n}[\mu_{1},\cdots,\mu_{n}]$, where
$l(\nu^{i})=n_{i}$ and $|\nu^{i}|=d_{i}$ for $i=1,2$.
\end{theorem}

\begin{remark}When we consider the simplest case $n=1$, the
above identity become the formula used in [6].
\end{remark}

\section{Some Examples of The Main Theorem}

In this section we give some examples of theorem 3.2.
\subsection{The case of $g=3$}If $g=3$, then $1\leq m\leq 3$. We
consider three cases.
\subsubsection{m=1}$LHS=-3\int_{\overline{\sM}_{3,1}}\lambda_{3}{\rm
ch}_{3}(\E)\psi_{1}$, and $3{\rm
ch}_{3}(\E)=\sum_{i+j=3}(-1)^{i-1}i\lambda_{i}\lambda_{j}=3\lambda_{3}-\lambda_{1}\lambda_{2},$
then we get
\begin{eqnarray*}
LHS&=&\int_{\overline{\sM}_{3,1}}\lambda_{3}(\lambda_{1}\lambda_{2}-3\lambda_{3})\psi_{1}=\int_{\overline{\sM}_{3,1}}\lambda_{1}\lambda_{2}\lambda_{3}\psi_{1}
,\\
RHS&=&b_{3}\frac{(-1)^{5}}{5}\binom {5}{0}\binom
{5}{4}B_{4}+\frac{1}{2}\sum_{g_{1}+g_{2}=3,g_{1},g_{2}>0}b_{g_{1}}b_{g_{2}}\frac{-1}{5}\binom
{2g_{2}-1}{0}\binom {5}{4}B_{4} \\
&=&-B_{4}(b_{3}+b_{1}b_{2}).
\end{eqnarray*}
Since $b_{1}=\frac{1}{24}, b_{2}=\frac{7}{5760},
b_{3}=\frac{31}{967680}$, we have
\begin{equation}
\int_{\overline{\sM}_{3,1}}\lambda_{1}\lambda_{2}\lambda_{3}\psi_{1}=\frac{1}{362880}.
\end{equation}
\subsubsection{m=2}In this case we have
$LHS=2\int_{\overline{\sM}_{3,1}}\lambda_{3}{\rm
ch}_{2}(\E)\psi^{2}_{1}$, and $2!{\rm
ch}_{2}(\E)=2\lambda_{2}-\lambda^{2}_{1}$, $B_{3}=0$. Then we have
\begin{eqnarray*}
LHS&=&\int_{\overline{\sM}_{3,1}}(2\lambda_{2}\lambda_{3}-\lambda_{3}\lambda^{2}_{1})\psi^{2}_{1}
,\\
RHS&=&b_{3}\sum^{1}_{k=0}\frac{(-1)^{4-k}}{5-k}\binom {5}{k}\binom
{5-k}{3}B_{3}
\\
&+&\frac{1}{2}\sum_{g_{1}+g_{2}=3,g_{1},g_{2}>0}b_{g_{1}}b_{g_{2}}\sum^{1}_{k=0}\frac{(-1)^{5-k}}{5-k}\binom
{2g_{2}-1}{k}\binom {5-k}{3}B_{3} \\
&=&0,
\end{eqnarray*}
hence
$$\int_{\overline{\sM}_{3,1}}\lambda_{3}\lambda^{2}_{1}\psi^{2}_{1}=2\int_{\overline{\sM}_{3,1}}\lambda_{2}\lambda_{3}\psi^{2}_{1}.$$
Using the formula
$\int_{\overline{\sM}_{3,1}}\lambda_{2}\lambda_{3}\psi^{2}_{1}=\frac{1}{120960}$,
we get
\begin{equation}
\int_{\overline{\sM}_{3,1}}\lambda_{3}\lambda^{2}_{1}\psi^{2}_{1}=\frac{1}{60480}.
\end{equation}
\subsubsection{m=3}In this case
$LHS=-\int_{\overline{\sM}_{3,1}}\lambda_{3}{\rm
ch}_{1}(\E)\psi^{3}_{1}$ and ${\rm ch}_{1}(\E)=\lambda_{1}$, hence
\begin{eqnarray*}
LHS&=&-\int_{\overline{\sM}_{3,1}}\lambda_{1}\lambda_{3}\psi^{3}_{1},
\\
RHS&=&b_{3}\sum^{2}_{k=0}\frac{(-1)^{5-k}}{5-k}\binom {5}{k}\binom
{5-k}{2}B_{2} \\
&+&\frac{1}{2}\sum_{g_{1}+g_{2}=3,g_{1},g_{2}>0}b_{g_{1}}b_{g_{2}}\sum^{{\rm
min}(2g_{2}-1,2)}_{k=0}\frac{(-1)^{2g_{2}-1-k}}{5-k}\binom
{2g_{2}-1}{k}\binom {5-k}{2}B_{2} \\
&=&-\frac{9}{2}b_{3}B_{2}-\frac{1}{2}b_{1}b_{2}B_{2} \\
&=&-\frac{41}{1451520},
\end{eqnarray*}
so
\begin{equation}
\int_{\overline{\sM}_{3,1}}\lambda_{1}\lambda_{3}\psi^{3}_{1}=\frac{41}{145120}.
\end{equation}

\begin{remark} The values of (37) and (39) match with the results in
[9], the identity (38) is a new result.
\end{remark}

\section{Acknowledgements}

The author would like to thank his advisor, Professor Kefeng Liu,
for his encouragement, helpful discussion on Hodge integrals and
correction of grammar mistakes. I also thank Professor Edna Cheung
who teach me a lot both in physics and in typing articles.
Professor Chiu-Chu Liu commented on this note and pointed out some
mistakes.

\end{document}